\documentclass[12pt]{article}

\usepackage{amssymb}

\newtheorem{satz}{Theorem}[section]
\newtheorem{defi}[satz]{Definition}

\newtheorem{kor}[satz]{Corollary}
\newtheorem{lem}[satz]{Lemma}

\newtheorem{bem}[satz]{Remark}

\newtheorem{prop}[satz]{Proposition}

\newtheorem{bsp}[satz]{Example}

\newcounter{Roma}
\newcounter{Ara}
\newcounter{let}

\newenvironment{Ar}{\begin{list}{(\arabic{Ara})\hfill}{\usecounter{Ara}
\labelsep3mm \leftmargin1.1cm \labelwidth8mm}}{\end{list}}

\newenvironment{Rom}{\begin{list}{(\roman{Roma})\hfill}{\usecounter{Roma}
\labelsep3mm \leftmargin1.1cm \labelwidth8mm}}{\end{list}}

\makeindex

\begin{document}

\newcommand{\nc}{\newcommand}

\nc{\mapco}{\,\colon\, }
\nc{\ab}{^{ab}}

\nc{\comment}[1]{
}

\nc{\catc}{{\C}}
\nc{\we}{\vee}

\nc{\hra}{\hookrightarrow}

\nc{\epi}{epimorphism}
\nc{\repi}{regular epimorphism}
\nc{\mono}{monomorphism}
\nc{\iso}{isomorphism}
\nc{\coker}[1]{\mbox{${\rm Coker}(#1)$}}
\nc{\Ker}[1]{\mbox{Ker$(#1)$}}
\nc{\defgl}{\stackrel{def}{=}}

\nc{\V}{\vspace{3mm}}
\nc{\VV}{\vspace{4mm}}
\nc{\lra}{\longrightarrow}
\nc{\lla}{\longleftarrow}
\nc{\mr}[1]{ \stackrel{#1}{\lra} }
\nc{\ml}[1]{ \stackrel{#1}{\lla} }
\nc{\hmr}[1]{\hspace{2mm}\stackrel{#1}{\lra}\hspace{2mm}}
\nc{\hml}[1]{\hspace{2mm} \stackrel{#1}{\lla}\hspace{2mm}}
\nc{\N}{\noindent}
\nc{\st}{^{\prime}}
\nc{\ot}{\otimes}
\nc{\hcong}{ \hspace{2mm}\cong\hspace{2mm}  }
\nc{\hfbox}{\hfill$\Box$}
\nc{\REF}[1]{(\ref{#1})}

\def\Z{\ifmmode{Z\hskip -4.8pt Z} \else{\hbox{$Z\hskip -4.8pt Z$}}\fi}

\def\Q{\ifmmode{Q\hskip-5.0pt\vrule height6.0pt depth
0pt\hskip6pt}\else{\hbox{$Q\hskip-5.0pt\vrule height6.0pt depth
0pt\hskip6pt$}}\fi}

\newcommand{\NN}{\mbox{$I\!\!N$}}

\nc{\Ph}{\phantom{}}

\nc{\BE}{\begin{equation}}
\nc{\EE}{\end{equation}}

\nc{\dst}{\displaystyle}
\nc{\sst}{\scriptscriptstyle}
\nc{\ssst}{\scriptscriptstyle}
\nc{\proof}{\N{\bf Proof\,:}\quad}
\nc{\proofof}[1]{\N{\bf Proof of {#1}\,:}\quad}
\nc{\proofofthm}[1]{\N{\bf Proof of theorem \ref{#1}\,:}\quad}

\nc{\htt}[1]{^{\otimes #1}}

\newcommand{\ssur}[2]{\mbox{$#1 \!\to\!\!\!\!\!\to\! #2$}}

\def\mapsel#1{\mbox{$\rule[-5mm]{0mm}{12mm} \searrow\rlap
{$\vcenter{\makebox[0mm][r]{$\scriptstyle#1$\hspace{12.2mm}}}$}$}}

\def\surltop#1{\makebox[1cm]{\mbox{$\stackrel{#1}{\makebox[0mm]{$\lla$}\hspace{0.7mm}
\makebox[0mm]{$\lla$}}$}}}

\nc{\Sur}[1]{\mbox{$\:\stackrel{#1}{\lra\!\!\!\!\!\to\,}\:$}}

\def\INJ{\mbox{\mathsurround=0pt
\makebox[0mm][r]{\parbox{0mm}{\rule[-0.65mm]{0mm}{0.2mm}$\scriptscriptstyle>$}}
\makebox[0.7cm][l]{\parbox{0.7cm} {$\lra$}}}}

\def\Inj#1{\mbox{$\:\stackrel{#1}{\INJ}\:$}}

\nc{\Injup}[1]{\mapup{#1}}
\nc{\injup}[1]{\mapup{#1}}
\nc{\injdown}[1]{\mapdown{#1}}


\def\mapup#1{\mbox{$\rule[-1mm]{0cm}{0.7cm}  
\makebox[0mm][r]{\raisebox{0.2mm}{$\scriptstyle\phantom{\cong}$}\hspace{0.6mm}}
\bigg\uparrow\rlap{$\vcenter{\hbox{$\scriptstyle#1$}}$}$}}

\def\isoup#1{\mbox{$\rule[-1mm]{0cm}{0.7cm}  
\makebox[0mm][r]{\raisebox{0.2mm}{$\scriptstyle\cong$}\hspace{0.6mm}}
\bigg\uparrow\rlap{$\vcenter{\hbox{$\scriptstyle#1$}}$}$}}

\def\mapdown#1{\mbox{$\rule[-1mm]{0cm}{0.7cm}  
\makebox[0mm][r]{\raisebox{0.2mm}{$\scriptstyle $}\hspace{0.6mm}}
\bigg\downarrow\rlap{$\vcenter{\hbox{$\scriptstyle#1$}}$}$}}

\def\isodown#1{\mbox{$\rule[-1mm]{0cm}{0.7cm}  
\makebox[0mm][r]{\raisebox{0.2mm}{$\scriptstyle\cong$}\hspace{0.6mm}}
\bigg\downarrow\rlap{$\vcenter{\hbox{$\scriptstyle#1$}}$}$}}

\def\isor#1{\mbox{$\smash{\mathop{\longrightarrow}\limits^{\cong}_{#1}}$}}

\def\isol#1{\mbox{$\smash{\mathop{\longleftarrow}\limits^{\cong}_{#1}}$}}

\def\surdown#1{\makebox[0mm]{$\mapdown{\raisebox{0.4mm}{$\scriptstyle#1$}}$}
\makebox[0mm]{\raisebox{-2.15mm}{$\downarrow$}}}

\def\surup#1{\makebox[0mm]{$\mapup{\raisebox{0.4mm}{$\scriptstyle#1$}}$}
\makebox[0mm]{\raisebox{0.7mm}{$\mapup{}$}}}

\newcommand{\surr}[2]{@\cdhgeneric>->->\twoheadrightarrow>#1>#2>}

\newcommand{\surl}[2]{@\cdhgeneric>\twoheadleftarrow>->->#1>#2>}

\newcommand{\brokenr}[2]{
@\cdhgeneric>\raise2.8pt\hbox to3.5pt{\hrulefill}\mkern9mu>
\raise2.8pt\hbox to3.5pt{\hrulefill}\hbox to5pt{}>
\mkern-7.5mu\dashrightarrow>#1 >
#2>}

\newcommand{\brokenup}[1]{
@\cdvgeneric>\hat\cdot
>\raisebox{3pt}{$\vdots$}>
\vbox{\kern3.5pt\hbox{$\cdot$}\kern-3.5pt}
> >\hspace{1mm}#1>}

\newcommand{\functlr}[2]{
\raisebox{0.4pt}{$\hss\begin{CD}
@>\vbox{\hbox to 0pt{$\hss\begin{CD}@<#1<<\end{CD}\hss$}\vskip-2pt}
>#2 >
\end{CD}\hss$}
}

\newcommand{\maprr}[2]{
\raisebox{-0.9pt}{$\hss\begin{CD}
@>\vbox{\hbox to 0pt{$\hss\begin{CD}@>#1>>\end{CD}\hss$}\vskip-3pt}
>#2 >
\end{CD}\hss$}
}

\newcommand{\mapdd}[2]{
@\cdvstandard>\downarrow\hspace{1.5pt}\kern-4pt\hspace{1.5pt}\Big\downarrow>#1>#2>}

\newcommand{\mapud}[2]{
@\cdvstandard>\uparrow\kern-3pt\Big\downarrow>#1>#2>}

\newcommand{\sepi}[3]{\,\mbox{$#1\,: \ssur{#2}{#3}$}\,}

\nc{\auf}{\twoheadrightarrow}
 
\nc{\ruled}{\rule[-4mm]{0mm}{0mm}}

\def\VC#1#2#3{\makebox[0cm][l]{\hspace*{-#1mm}\makebox[#2cm]{\raisebox{4mm}{$#3$}}}
\begin{minipage}{0cm} \unitlength1cm \hspace*{-#1mm}
\begin{picture}(1,0) \put(0,0){\vector(1,0){#2}}\end{picture}  \end{minipage}
}

\nc{\qu}{quadratic}
\nc{\GBoGB}{G/BG\st \otimes G/BG\st}
\nc{\lstar}{_{\raisebox{-1mm}{$*$}}}
\nc{\sm}{\:{ \wedge}\:}
\nc{\hsm}[1]{^{\sst \wedge #1}}

\nc{\Rmod}{${\bf R}$-module}

\nc{\map}[3]{\mbox{$#1 \mapco #2 \to #3$}}

\nc{\rond}{{\,\sst \circ\,}}

\nc{\ruleu}{\rule{0mm}{7mm}}

\nc{\T}[1]{\tilde{#1}}

\nc{\Imm}[1]{\mbox{${\rm Im}(#1)$}}

\nc{\tw}{\end{document}}

\newcommand{\Nat}[1]{\mbox{$I\!\!N^{#1}$}}

\nc{\dsg}{dimension subgroup}

\nc{\Pz}{\mbox{$P_{2,\Z}^{\cal N} $}}
\nc{\kpz}{\mbox{$p_{2,\Z}^{\cal N} $}}

\nc{\Pd}{\mbox{$P_3^{\cal N} $}}
\nc{\kpd}{\mbox{$p_3  $}}

\nc{\Pn}{\mbox{$P_{n,R}^{\cal N} $}}
\nc{\kpn}{\mbox{$p_{n,R}^{\cal N}  $}}

\nc{\Dn}{\mbox{$D_{n,R}^{\cal N} $}}
\nc{\Dz}{\mbox{$D_{2,R}^{\cal N} $}}
\nc{\Dd}{\mbox{$D_{3,R}^{\cal N} $}}

\nc{\In}{\mbox{$I^n_{R,\cal N}  $}}
\nc{\kpzr}{\mbox{$p_{2,R}  $}}

\vspace*{0.7cm}

\begin{center} \Large \sc    The relative second Fox and third dimension subgroup of arbitrary groups \normalsize  
 \vspace{0.8cm}

\renewcommand{\thefootnote}{ }
\normalsize \bf  \sc  Manfred Hartl\normalsize \footnote[1]{Research partially supported by an individual
fellowship of the Human Capital and
Mobility }
\footnote[1]{Programme of the European Union during the years 1993-1995.}
\end{center}

\vspace{8mm}

\begin{abstract}

\N Let $I_R(G)$ denote the augmentation ideal of the group algebra $R(G)$ of a group $G$ with coefficients in a commutative ring
$R$. We give a complete description of the third relative dimension subgroup $G\cap(1+I_R(K)I_R(G)+I^3_R(G))$ and the second
relative Fox subgroup $G\cap(1+I_R(K)I_R(H)+I^2_R(G)I_R(H))$ for any subgroups $K$ and $H$ of $G$.

\end{abstract}\vspace{3mm}

\N Mathematical subjects classification:  20C07 (Primary), 20J05 (Secondary).
\vspace{4mm}

\section{Introduction}

It is a classical problem to study the link between filtrations of groups and of their group algebras. There are
mainly two cases much studied in the literature, namely dimension subgroups and Fox subgroups. As to the first, 
fix   a commutative ring of coefficients,
$R$,  with unit $1_R$. Now an appropriate filtration of a group $G$ (by an N-series ${\cal N} =
\{N_i\}$ or a restricted N-series) gives rise to a descending filtration of the group
algebra $R(G)$    by ideals $\In(G)$, see \cite[III.1.5]{Pa} or section 2 below. For example,  the lower central series
$\gamma\,\colon\, G=G_1\supset G_2\supset\cdots $ of $G$ induces the filtration by the powers $I_R^n(G)$ of 
the augmentation ideal $I_R(G)$ of $R(G)$. Pulling back the induced filtration $\{\In(G)\}$ to $G$
 defines the series of {\em dimension subgroups with respect to $R,\cal N$}\/, $\Dn (G) =  G \cap
(1 + I_{\cal N}^n(G))$. Now  the so-called {\em dimension subgroup problem}\/ asks 
 whether the canonical inclusion $N_n
\subset \Dn (G)$ is an equality. For classical \dsg s
$D_n^{\gamma} (G)$ this is true for $n\le 3$ and is now known to be false for all $n\ge 4$, by means
of counterexamples due to Rips \cite{Rips} ($n=4$)and later by Gupta \cite{GuDn} ($n\ge 4$).  While their method is
combinatorial, a {\em homological approach}\/ to the problem was inaugurated by Passi, from where the
notion of \dsg s {\em relative to a subgroup}\/ $K$ of $G$ emerged in a natural way: these are
defined by $\Dn (G,K) =  G \cap (1 + I_R(K)I_R(G)+ I_{R,\cal N}^n(G))$. Since then, they proved more
and more to be an appropriate tool  for the study  of the classical case $K=\{1\}$. Many qualitative
properties of relative dimension subgroups were established,  notably by  
Kuz'min \cite{Ku}. Nevertheless, the main problem, concerning their explicit computation, remains difficult,
even in low dimensions $n$. With regard to the ring of coefficients,  only $D_{n,R}^{\gamma}(G)$,
$n\le 3$, seems to be known for {\em all}\/ commutative rings $R$, by work of Sandling \cite{Sa}. Focussing on $R=\Z$,  the
``relative" analogue of the classical \dsg\ problem asks when the  inclusion $K_2N_n \subset
D_{n,\Z}^{\cal N}(G,K)$ is an equality. This is true for $n\le 2$, and various conditions on $G$ and $K$ were
exhibited in the literature ensuring that it also holds for $n=3$. In the first part of this paper we
compute $\Dd (G,K)$ in general and show that equality does {\em not}\/ always hold,
\comment{ (even for ${\cal N} =  \gamma$),}
 the minimal counterexample for ${\cal N}=\gamma$ and $R=\Z$
being of order $2^6$. We note that the computation of $D_{3,Z\!\!\!Z}(G,K)$, achieved and distributed by the author in 1994, was
reproved by Tahara, Vermani and Razdan by different methods in \cite{Razdan}.

In the second part of the paper we unify the study of relative \dsg s with the one of {\em Fox
subgroups}\/ $G \cap (1_R + I^n_R(G)\, I_R(H))$ for a subgroup $H$ of $G$. 
After a long history (for a review see \cite{Gu} or \cite{Gu-Cu}), they are now completely known for
free groups $G$ (${\cal N} =\gamma$, $R= \Z$), thanks to work of N.\ Gupta, Hurley and Yunus.
 The case of arbitrary groups is much harder;   for
$n=2$ (and $R=\Z$) the problem was  solved by K.\ Gupta and  M.\ Curzio in \cite{Gu-Cu} but seems to be completely open for $n>2$.
On the other hand, a first case of a relative version of Fox subgroups was considered in \cite{Razdan} where the group $G\cap
(1+ I_{Z\!\!\!Z}(K)I_{Z\!\!\!Z}(H)+I_{ Z\!\!\!Z }^2(G)I_{Z\!\!\!Z}(H) )$ is determined for $H$ normal and $K$ being a specific
subgroup of $G$ containing $[H,G]$. Generalizing  the two
last-mentioned  results  we here determine the group $G \cap (1_R +  I_R(K)\,I_R(H) +
I^2_R(G)\,I_R(H))$ for any subgroups $H,K$ of $G$ and coefficient rings $R$.

Although this result formally includes the case $H=G$ treated in section 2 by a different method, the
result obtained there is much simpler and does not seem to be easily deducible from the  more general
formula in section 3. Nor it seems possible to generalise the method of section 2 to the latter
case, because of the essential difference of the behaviour of additive and of quadratic functors
(like $\bigotimes^2 , SP^2$) with respect to subgroups.

\comment{
Dimension subgroups of a group $G$ {\em relative to a subgroup}\/ $H \le G$ are defined by
$\Dn (G,H) = G \cap (1 + I(H)I(G)+ I_{\cal N}^n(G))$, where the intersection is taken in the group ring
$\Z(G)$ with augmentation ideal $I(G)$. This generalisation of classical dimension subgroups
$\Dn (G)=\Dn (G,\{1\}) $ emerged from the work of Passi, in relating the groups $\Dn (G)$ to a
`polynomial approximation' of the Schur multiplier $H_2(G)$ and thus to homological algebra. Since
then, many qualitative properties of relative dimension subgroups were established,  
notably by Gupta and Kuz'min. Nevertheless, the main problem, concerning their explicit
determination, remains difficult, even in low dimensions $n$. In fact, the `relative' analogue of the
classical \dsg\ problem asks when the canonical subgroup $H_2N_n \subset \Dn (G,H)$ is an equality,
where $G=N_1\supset N_2 \supset\cdots$ denotes the lower central series. For classical \dsg s
$\Dn (G)$ this is true for $n\le 3$ and is now known to be false for all $n\ge 4$, by means of
counterexamples due to Rips and Gupta. In the relative case, the relation $H_2N_n \subset \Dn (G,H)$ is
true for $n\le 2$, and various conditions on $G$ and $H$ were exhibited in the literature ensuring
that this also holds for $n=3$. In this paper we determine $\Dd (G,H)$ in general and show that
equality does {\em not}\/ hold in general, the minimal counterexample for ${\cal N}=\gamma$ being of
order $2^6$. Indeed, the default is made precise by the following theorem.\V

Our method is based on two rather elementary homological observations which seems to generalise to
$n=4$, based on the results of \cite{Q3}.
}

\section{The third relative dimension subgroup}

Recall the notation from the introduction. In particular,  recall that an {\em N-series}\/ ${\cal N} = \{N_i\}$ is a descending
chain of subgroups
  \[  G = N_1 \supset N_2 \supset \ldots \supset 1 \]
such that $[N_i,N_j] \subset N_{i+j}$ (with $[a,b] = aba^{-1}b^{-1}$ for $a,b \in  G$). 
An N-series ${\cal N}$ induces a descending chain of two-sided ideals 
  \[  R(G) \supset I^1_{R,\cal N} \supset I^2_{R,\cal N} \supset \ldots \supset 0 \]
by defining $I^n_{R,\cal N}$ to be the $R$-submodule of $R(G)$ generated by the elements
  \[  (a_1 -1) \cdots (a_r -1)\,, \quad \mbox{$a_i \in N_{k_i}$, such that $k_1+ \ldots +
k_r  \ge n$.}  \] 
 For a subgroup $K$ of $G$ and $n\ge 1$ define
  \[ \Dn (G,K) =  G \cap (1 + I_R(K)I_R(G)+ I_{R,\cal N}^n(G))\:. \]
Note that for $n=3$  the  case of an arbitrary N-series ${\cal N}$ can be reduced to the case  ${\cal N}=
\gamma $ by the identity $\Dd (G,K)/N_3 = D_3^{\gamma}(G/N_3,\,K N_2/N_3)$, but we do not make use of
this reduction as our method  genuinly works for an arbitrary N-series. 

The main goal   of this section is to prove the following result.
\vspace{2mm}

\begin{satz}\label{D3R}\quad Let G be a group, K a subgroup, ${\cal N}$ an N-series of G and  R a
commutative ring with unit $1_R$. Then
  \[ D_{3,R}^{\cal N}(G,K)  = U_0  N_3   Z_2   \prod_{\matrix{\scriptstyle  p \in  \sigma(R) \cr \scriptstyle p\:
odd}}  t_p(  \mbox{$G$ {\rm mod} $ U_0 N_3 $})   \cap (U_{p^e} N_3 G^{p^e})  \]
where
\begin{itemize}

\item  $U_m= sgp\{ [a,b^k]\,\|\, a,b \in G,  k\in \Z, \; a^k,b^k \in KN_2 G^m\}$;

\item  $\sigma(R) = \{p|  \mbox{$p$ is a prime and $ p^nR$} =p^{n+1}R$ for some
$n\ge 0\}$, and for $p \in \sigma(R)$, $p^{e}  = p^{e(p)}$ is the smallest power of $p$ for
which $p^{e }R = p^{e+1}R$;

\item  $ t_p(  \mbox{$G$ {\rm mod} $ U_0 N_3$}) = \{g \in G\,|\,\mbox{$g^{p^k} \in U_0 N_3$
for some $k\ge 0$}\}$;
   
\item $Z_2 =\{1\}$ if $2\notin \sigma(R)$, else 
$Z_2 =  t_2(  \mbox{$G$ {\rm mod} $ U_0 N_3 $})   \cap (U_{2^{e(2)}} N_3 G^{2^{e(2)+1}} 
V^{2^{e(2)}} )  $, where  $V=\{g\in G\,\|\, g^{e(2)-1} \in KN_2 G^{2^{e(2) }}\} $   if $e(2) >0$
and $V=G$ else.

\end{itemize}
\end{satz}
\vspace{2mm}

For the proof one first reduces to the case $R=\Z/m\Z$ via 
the "universal coefficient decomposition" obtained
in theorem 1 of \cite{Unico}. Then what remains to prove is the following.

\begin{satz}\label{D3}\quad Let $G$ be a group, $K$ a subgroup, ${\cal N}$ an N-series of G and $m\ge
0$.  If $m$ is even let $V=\{a\in G\,\|\,a^{m/2} \in KN_2G^m\}$. Then  
\[  D_{3,\Z}^{\cal N}(G,K) = U_0 N_3 \] 
\[  D_{3,\Z/m\Z}^{\cal N}(G,K) = U_mN_3G^m  \quad \mbox{if m is odd} \] 
\[ \makebox[14.7cm]{ \phantom{$\Box$} \hfill $ D_{3,\Z/m\Z}^{\cal N}(G,K) =
U_mN_3G^{2m}V^m  \quad \mbox{if m is even.} $\hfbox } \] 
\end{satz}
For $K=\{1\}$ one rediscovers a result of
Sandling \cite{Sa}. Some other special cases are resembled in the following

\begin{kor}\label{D3=}\quad
The inclusion $K_2N_3 \subset D_{3,\Z}^{\gamma}(G,K)$ is an {\bf equality} if one of the following
conditions holds:

\begin{Ar}

\item $[K,G] \subset G_3$ (cf.\ \cite[Theorem 5.9]{Pa});

\item there exists a normal subgroup $N\subset G$ such that $G=NK$ and $N\cap K$ is central in $G$
(cf.\ \cite{Ka-Ve});

\item $K$ is normal and $G/K$ is cyclic (cf.\ \cite[Theorem V.5.4]{Pa}); 

\item  one of the following three groups is torsion-free: $G/K_2G_3$ (cf.\cite{Ka-Ve}), \\
$[G,K]G_3/K_2G_3$ or $G/KG_2$;

\item  the abelian group $KG_2/G_2$ is
divisible.\hfbox

\end{Ar}
\end{kor}

These facts are easily  derived  from theorem  \ref{D3} or the -- essentially equivalent -- theorem
\ref{P2tau}.

 Nevertheless, theorem \ref{D3} surprisingly shows that the inclusion 
$K_2G_3 \subset D_{3,\Z}^{\gamma}($\linebreak$ G, K)$ is   {\bf not  always} an equality, as was
suggested by the known partial results reviewed in \ref{D3=}. Indeed, we find counterexamples which
are $p$-groups for {\bf any prime} $p$, see \ref{cex} below; this is in contrast to the case of
classical dimension subgroups (i.e.\ $K=1$) which coincide with the terms of the lower central
series of $G$ unless $p=2$ (due to a recent result of N.\ Gupta).

\begin{bsp}\label{cex}\quad \rm Let $p$ be a prime and $0<r\le s$.
Define
  \[ G=\langle x,y\, |\,1= x^{p^{s+1}} = y^{p^{s+1}} = [x,[x,y]] = [y,[x,y]]\,\rangle \:.\]
Let $K = sgp\{x^{p^r}, y^{p^{s}}, [x,y]\}$. Then $z = [x,y]^{p^{s}} =
[x,\,y^{p^{s}}] \in D_{3,\Z}^{\gamma}(G,K)$ by \ref{D3}, but $z$ has order $p$ modulo
$K_2G_3 = \{1\}$.\hfbox

\end{bsp}

In order to prove theorem \ref{D3}, we need to study a related quotient of the group  
algebra. The {\bf relative polynomial group}\/ is defined by
  \[ \Pn  (G,K) = I_R(G)/(I_R(K)I_R(G)+ I_{R,\cal N}^{n+1}(G)) \,.\rule[-3mm]{0mm}{3mm} \]
It generalizes the well-known  polynomial
group of Passi, $\Pn (G) = P_{n,R}^{\cal N} (G,$\linebreak $\{1\})$, see \cite{Pa68} or also
\cite{Pa}. Also the relative version is implicit in the work of Passi and  various other  places in
the literature. For a discussion and more properties, in particular of the torsion subgroup and the
torsion-free quotient of  $P_{n,\Z}^{\gamma}(G,K)$, see \cite{PolProp}. Here we only resemble some
elementary properties in the following

\begin{lem}\label{Pnsequ}\quad Let $K\lhd G$ be a normal subgroup. Then $\Pn (G,K)$ admits a left
$R(G/K)$-module structure induced by multiplication in $ I_R(G)$, which makes the canonical map
  \[ \map{p_{n,R}^{\cal N}}{G}{\Pn (G,K)}\,,\quad a \mapsto (a-1)+ I_R(K)I_R(G)+ I_{R,\cal
N}^{n+1}(G)  \]
into a (left) derivation, $\kpn (ab) = a\kpn (b) + \kpn (a)$ for $a,b \in G$. Moreover,
there is an exact sequence of  $R(G/K)$-linear homomorphisms
  \[ R\ot (KN_{n+1}/K_2N_{n+1}) \mr{p_{n,R}^{\cal N} i} \Pn (G,K) \:\mr{P_{n,R}^{\cal N}(\pi)}\:
P_{n,R}^{\pi\cal N} (G/K) \,\to\,  0\,,\] 
where $i$ is induced by the inclusion $K \hra G$, $p_{n,R}^{\cal N} i$ also
denotes its $R$-linear extension, ${\pi\cal N}$ is the N-series ${\pi\cal N}_i = \pi(N_i)$ with
$\sepi{\pi}{G}{G/K}$, $P_{n,R} ^{\cal
N}(\pi)\{a-1\} = \{aK -1\}$, and where the $R(G/K)$-action on the left-hand term is induced by
$R$-linear extension of the $(G/K)$-action induced by conjugation in $G$. \hfill $\Box$ \end{lem}

Note that by definition,
  \[  \Ker{p_{n,\Z}^{\cal N} i} = \frac{\dst KN_{n+1} \cap D_{n+1,\Z}^{\cal N}(G,K)}{\dst K_2
N_{n+1}} \] which is nontrivial in general; indeed, in the case $K=N_n$ one rediscovers the
classical dimension subgroup problem, so that in taking $G$ to be the group of Rips \cite{Rips},
$K=G_3$ and ${\cal N} = \gamma$ we have  $\Ker{p_{3,\Z}^{\gamma}i} \neq \{1\}$. In the next theorem  we calculate
$\Ker{p_{2,\Z}^{\cal N} i}$ by extending the sequence in Lemma \ref{Pnsequ} on the left, as follows. 

Consider the following part of a six-term exact sequence for the tensor and torsion product of
abelian groups.
\comment{
\\
\N\makebox[14.7cm]{ \makebox[0mm]{
\begin{minipage}{20cm}\small
\[ {\rm Tor}_1^{\Z}(G/KN_2, G/KN_2) \mr{\tau}   (G/KN_2) \ot (KN_2/N_2)   \mr{id \ot j} (G/KN_2) \ot
(G/N_2) \mr{id \ot q}  (G/KN_2)\ot (G/KN_2) \,\to\,0 \]
\end{minipage}\ruled
}\rule[-5mm]{0mm}{3mm} }
\rule{0mm}{7mm}
}
\comment{
  \[ {\rm Tor}_1^{\Z}(G/KN_2, G/KN_2) \mr{\tau}   G/KN_2 \ot KN_2/N_2   \mr{id \ot j} G/KN_2\ot G/N_2
\mr{id \ot q}  G/KN_2\ot G/KN_2 \,\to\,0 \]
}
\\
\N\makebox[14.7cm]{ \makebox[0mm]{
\begin{minipage}{20cm}\small
\[ \begin{array}{ccrcc}
{\rm Tor}_1^{\Z}(G/KN_2, G/KN_2) &\mr{\tau} &  (G/KN_2) \ot (KN_2/N_2)   &\mr{id \ot j} &(G/KN_2)
\ot (G/N_2) \ruled\\
& & & & \mapdown{id \ot q} \\
& & 0 & \leftarrow & (G/KN_2)\ot (G/KN_2) \rule{0mm}{6mm}
\end{array} \]
\end{minipage}\ruled
}\rule[-5mm]{0mm}{3mm} }
\rule{0mm}{7mm}
 Here $j,q$ are the canonical inclusion and
quotient map, respectively. Moreover, commutation in $G$ induces a homomorphism
 \[ \map{[\,,\,]}{ (G/KN_2) \ot (KN_2/N_2)}{KN_3/K_2N_3}\,.\]
\vspace{2mm}

\begin{satz}\label{P2tau}\quad Let $K$ be a normal subgroup of a group $G$. Then the following
sequence of natural homomorphisms is exact:
  \[ {\rm Tor}_1^{\Z}(G/KN_2, G/KN_2) \mr{[\,,\,] \tau} KN_{3}/K_2N_{3} \mr{p_{2,\Z}^{\cal N} i}
\Pz (G,K) \mr{P_{2,\Z}^{\cal N}(\pi)} P_{2,\Z}^{\pi\cal N} (G/K) \to 0\,.\]
\end{satz}

We remark that this result admits an application in group cohomology with respect to the variety
of $2$-step nilpotent groups, thus solving a problem of Leedham--Green. This will be presented
elsewhere in a more general context.

The proof of theorem \ref{P2tau} requires a homological lemma which is useful also elsewhere.

\begin{lem}\label{cruclem}\quad Let A be an abelian group and $B\stackrel{j}{\hra}A$ a subgroup.
Consider the following homomorphisms
  \[ (A/B) \ot B \surltop{q\ot id} A \ot B  \mr{\nu} A \sm A \mr{\ell} A\ot A \Sur{q\ot id} (A/B) \ot A
\] where $A\sm A = A\ot A/sgp\{a\ot a \| a \in A\}$, $a\sm a\st = \{a \ot a\st\}$, $\nu(a\ot b) = a
\sm b$, $l(a\ot a\st) = a\ot a\st - a\st \ot a$, and where $q$ is the quotient map. Then
  \[ \Ker{(q\ot id)\ell} = \nu(q\ot id )^{-1} \Imm{\map{\tau}{{\rm Tor}_1^{\Z}(A/B,A/B)}{(A/B)\ot
B}}\,,\] where $\tau$ appears in the following part of a six-term exact sequence,
\[  {\rm Tor}_1^{\Z}(A/B,A/B) \mr{\tau} (A/B) \ot B \mr{id \ot j} (A/B) \ot A \mr{id \ot q} (A/B)
\ot (A/B) \to 0 \,.\]
\end{lem}

\proof  Consider the following commutative square
  \[ \matrix{
A \sm A & \mr{\ell}  & A \ot A  & \mr{q\ot id}  &  A/B \ot A  \rule[-3mm]{0mm}{3mm} \cr
\surdown{q\ot q}  &  &  &  & \surdown{id \ot q}  \cr
(A/B) \sm (A/B)  & \VC{5}{3.5}{\ell} & & & (A/B) \ot (A/B)
}\]
As the map $\ell$ is injective for all abelian groups (see \cite{Ba-GrI}), $\Ker{(q\ot id)\ell}$ is
contained in $\Ker{q\ot q} = \Imm{\nu}$, whence 

  \begin{eqnarray*}
\Ker{(q\ot id)\ell} &=& \nu\,\Ker{(q\ot id)\ell\nu} \\
  &=&  \nu\,\Ker{q\ot j} \\
  &=& \nu\,\Ker{(id\ot j)(q\ot id)}   \\
  &=& \nu(q\ot id)^{-1}\Ker{id \ot j}  \\
  &=&  \nu(q\ot id)^{-1} \Imm{\tau}\,.
\end{eqnarray*}
\Ph\hfill$\Box$\V

\proofofthm{P2tau}  Consider the following commutative diagram of homomorphisms
  \[\matrix{
(G/N_2) \sm (G/N_2)  &  \mr{\ell}  &  (G/N_2) \ot (G/N_2) & \Sur{q\ot id}  & (G/KN_2) \ot (G/N_2)
\rule[-3mm]{0mm}{3mm} \cr
\mapdown{c}  &  & \mapdown{\mu}  &  & \mapdown{\bar{\mu}}  \cr
N_2/N_3 & \mr{p_{2,\Z}^{\cal N}} &I_{\Z, \cal N}^2(G)/I_{\Z, \cal N}^3(G)  &  \Sur{\bar{q}}  & 
\frac{\dst \ruleu I_{\Z, \cal N}^2(G)  }{\dst  I_{\Z} (K)I_{\Z} (G) + I_{\Z, \cal N}^3(G)} 
}\]
where 
 \[ c(aN_2 \sm bN_2) = [a,b]N_3\,,\quad \mu(aN_2 \ot bN_2) = (a-1)(b-1) + I_{\cal N}^3(G)\,, \]
and where $\bar{q}$ is the canonical quotient map and $\bar{\mu}$ is induced by $\mu$. Indeed,
$\Ker{\bar{q}} = \mu\,\Ker{q \ot id}$, so the right-hand square is a pushout of abelian groups
(cf.\ \cite{Rot}), as is the left-hand square by the identity $Q_{2,\Z}^{\cal N}(G) = {\rm U}_2{\rm L}
^{\cal N}(G) =
(G/N_2) \ot (G/N_2)/l\,\Ker{c}$  obtained in \cite{Q3}, see also \cite{Ba-GrII}  (and which can also
be derived from Passi's theorem that $D_3^{\gamma}(G,\zeta_1(G))=G_3$, see \cite[V.5.9]{Pa}, by using
his technique in \cite[VIII.8.7]{Pa}). So by general nonsense (gluing of pushouts, which is easily
verified by using the universal property), also the exterior rectangle is a pushout. Therefore, 
  \[ (G
\cap (1 +I_{\Z}(K) I_{\Z}(G) + I_{\Z,\cal N}^3(G) ))/N_3 =  \Ker{\bar{q}\kpz} = c \,\Ker{(q\ot id)\ell}
\] 
where the first identity follows from the elementary relations  
  \BE\label{absch}  G \cap (1
+I_{\Z}(N)I_{\Z}(G) + I_{\Z,\cal N}^3(G) ) \:\subset\:  G \cap (1 +  I_{\Z,\cal N}^2(G))  = N_2\,. 
\EE 
 Now apply Lemma \ref{cruclem} for $A=G/N_2$ and $B=KN_2/N_2$. Then one has 
 \[ c\nu\,\Ker{q\ot id} = c\nu((KN_2/N_2) \ot (KN_2/N_2)) =  (K_2N_3)/N_3 \,,\]
so the result follows from commutativity of the following diagram.
    \[ \matrix{
(G/KN_2)\ot (KN_2/N_2)  &  \surltop{q\ot id}  &  (G/N_2) \ot (KN_2/N_2)  & \mr{\nu}  & (G/N_2) \sm
(G/N_2) \rule[-3mm]{0mm}{3mm} \cr 
\mapdown{[\,,\,]}  & \searrow &  &  &  \mapdown{c}  \cr
((KN_3)\cap N_2)/K_2N_3  &  \hra  & N_2/K_2N_3  & \surltop{}  &  N_2/N_3 \rule{0mm}{7mm}
}\]
\Ph \hfill $\Box$\V

\proofofthm{D3}  We first observe that $K$ may be replaced by the normal subgroup $KN_2G^m$, indeed,
  \[   I_R(KN_2G^m) \,\equiv  \, I_R(K) + I_R(N_2)  +  m\,I_R(G) \,\equiv  \, I_R(K) \quad \mbox{ 
mod $I^2_{R,\cal N}(G)$}\,,   \]
whence 
  \BE\label{Krepl} D_{3,R}^{\cal N}(G,K) = D_{3,R}^{\cal N}(G,KN_2G^m) = N_2G^m \cap (1_R + 
I_R(KN_2G^m) I_R(G)  + I_{R,\cal N}^3(G) )  \EE
where the second equation is obtained from the estimate
  \[ G \cap (1_R + I_R(KN_2G^m) I_R(G)  + I_{R,\cal N}^3(G)) \subset 
 G \cap (1_R +  I_{R,\cal N}^2(G))  =  G \cap (1_R + I_R^2(G)  + I_{R}(N_2))  \]
together with the isomorphisms
   \[  I_R(G)/( I_R^2(G) + I_R(N_2)) \,\cong\, R \ot (G/N_2) \,\cong\,  G/N_2G^m\,,\]
cf.\ \REF{n=2} and \REF{XotR} below.
 Together with  \REF{absch} and theorem \ref{P2tau} we obtain identities
  \begin{eqnarray}
 \frac{\dst D_{3,\Z}^{\cal N} (G,KN_2G^m)}{\dst K_2 N_3 (G^m)_2}  &=&  \frac{\dst
KN_2G^m \cap (1 +  I_{\Z}(KN_2G^m)I_{\Z}(G) + I_{\Z,\cal N}^3(G))}{\dst K_2N_3(G^m)_2} \nonumber\\
  &=&
\Ker{p_2^{\cal N} i}  \nonumber\\
  &=& {\rm Im}\Big( [\,,\,]\tau\,\colon\,{\rm Tor}_1^{\Z}(G/(KN_2G^m),
G/(KN_2G^m)) \mr{}  \nonumber\\
  & & \hspace{4cm} (KN_2G^m)/(K_2N_3(G^m)_2) \Big)  \label{D3=tau}
\end{eqnarray}
  To make the last term explicit we use the description
of the torsion product of abelian groups and the connecting homomorphism given in \cite[V.6]{ML}. In
fact, let $\langle \bar{x}_1,k, \bar{x}_2\rangle$ be a generator of Tor$_1^{\Z}(G/KN_2G^m,G/KN_2G^m)$, 
i.e., $k\in \Z$, $\bar{x}_i = x_i KN_2G^m$ for $x_i \in G$ such that $x_i^k \in KN_2G^m$,
$i=1,2$.  Then $[\,,\,]\tau \langle \bar{x}_1,k, \bar{x}_2\rangle = [x_1,x_2^k]$, so 
  \BE\label{tau=U} \Imm{ [\,,\,]\tau} =  U_mN_3 / K_2 N_3 (G^m)_2\,,\EE
noting that $U_m$ contains $K_2 (G^m)_2$. By \REF{Krepl} and \REF{D3=tau} this proves the assertion
for $m=0$.  \V

\N{\bf The case $m>0$.}\/ In the sequel we shall frequently use the canonical
identifications
  \BE\label{XotR} R \ot X \:\cong\: X/mX \:\cong\:  X\ot R\,,\quad  1_R \ot x \mapsto x + mX \mapsto
x\ot 1_R \EE
for abelian groups $X$, and the fact that ${m \choose 2} (x\ot 1_R) =   \mbox{} - 
{m \choose 2} (x\ot 1_R)$.\V

We start with the observation that $P_{2,R}^{\cal N}(G,KN_2G^m)\,\cong \, R \ot P_{2,\Z}^{\cal
N}(G,KN_2G^m)$. Now consider the following commutative diagram.
\V

\N\makebox[14.7cm]{ \makebox[0mm]{
\begin{minipage}{20cm}\small
  \[ \matrix{
{\rm Tor}_1^{\Z}(R,{\rm SP}^2(G/KN_2G^m))  &  \mr{\tau_1}  &  R \ot 
\frac{\dst  (G/N_2) \sm (G/N_2)}{\dst \Ker{(q\ot id)\ell}\ruled}  &  \mr{R \ot (q\ot id)\ell}  &  R \ot
(G/KN_2G^m)  \ot (G/  N_2)  \cr
\mapdown{\mu_{\ast}}  &  & \mapdown{R \ot c}  &  &  \mapdown{R \ot \bar{\mu}} \cr
\rule[-6mm]{0mm}{3mm}\ruleu {\rm Tor}_1^{\Z}(R, P_{2,\Z}^{\gamma} (G/KN_2G^m))  &  \mr{\tau_2}  &  R
\ot (KN_2G^m/U_mN_3)  & \mr{R \ot p_{2,\Z}^{\cal N}i}  &  R \ot P_{2,\Z}^{\cal N}(G, KN_2G^m) \cr
\| &  &  & & \cr
\ruleu {\rm Tor}_1^{\Z}(R, P_{2,\Z}^{\gamma} (G/KN_2G^m))  & \mr{\rho_{\ast}}  & 
{\rm Tor}_1^{\Z}(R,  G/KN_2G^m )  &    \mr{\tau_3}  &  R \ot {\rm SP}^2(G/KN_2G^m)
}\]
\end{minipage}\ruled
}\rule[-5mm]{0mm}{3mm} }
\rule{0mm}{7mm}

The lines are parts of the six-term exact sequences associated with the short exact sequences
  \[ ((G/N_2) \sm (G/N_2))/\Ker{(q\ot id)\ell}  \Inj{(q\ot id)\ell} (G/KN_2G^m) \ot (G/N_2) \Sur{q}
{\rm SP}^2(G/KN_2G^m)  \]
   \[ (KN_2G^m)/(U_mN_3)  \Inj{p_{2,\Z}^{\cal N}i} P_{2,\Z}^{\cal N}(G, KN_2G^m) 
\Sur{P_{2,\Z}^{\cal N}(\pi)}
P_{2,\Z}^{\gamma}(G/ KN_2G^m)  \]
  \BE\label{murho} \ruled\ruleu {\rm SP}^2(G/KN_2G^m)  \Inj{\mu} P_{2,\Z}^{\gamma}(G/ KN_2G^m) 
\Sur{\rho} G/ KN_2G^m  \EE 
In the first sequence, $q$ denotes the canonical quotient map.
The second sequence is obtained from \ref{Pnsequ}, \REF{D3=tau} and \REF{tau=U}. The third sequence
is due to Passi \cite{Pa69}, where SP$^2$ denotes the symmetric  tensor product, and
$\rho p_{2,\Z}^{\gamma}(gKN_2G^m) =  gKN_2G^m$.

Our goal is   to compute $\Ker{R\ot p_{2,\Z}^{\cal N} i} = \Imm{\tau_2}$. First note that
 tensoring with $R = \Z/m\Z$ leaves the first sequence unchanged  by \REF{XotR}; this implies that
$\tau_1=0$. So $\tau_2$  factors through a map
  \[ \map{\bar{\tau}_2}{ \Ker{\tau_3} = \Imm{\rho_{\ast}}\,\cong\,  {\rm coker}(\mu_{\ast})}{R \ot 
(KN_2G^m / U_m N_3)} \,, \]
so that 
 \BE\label{Imtau} \Imm{\tau_2} = \Imm{\bar{\tau}_2}\,.  \EE
In order to calculate \Ker{\tau_3} we use the canonical identification of Tor$_1^{\Z}(\Z/m\Z\,,A)$
with the subgroup $A_{(m)}$ of $m$-torsion elements of an abelian group $A$.\vspace{2mm}

\begin{lem}\label{Kertau3}\quad For an abelian group $A$ and $m\ge 0$, one has 
  \[ \Ker{\map{\tau_3}{A_{(m)}}{(\Z/m\Z) \ot {\rm SP}^2(A)}} \:=\:\left\{
\begin{array}{ll}
  A_{(m)} & \mbox{if $m$ is odd;} \\
( A_{(m)} \cap A^2) A_{(m/2)}  & \mbox{if $m$ is even.}
\end{array}
\right.  \]
\end{lem}

\proof The assertion is true for $m=0$ by exactness of sequence \REF{murho} and since for $m=0$,
$A_{(m)} = A = A_{(\frac{m}{2})}$. So suppose $m >0$. For $a \in A_{(m)}$ one has 
  \begin{eqnarray*}
\tau_3(a)  &=&  1_R \ot \mu^{-1}(m p_{2,\Z}^{\gamma}(a))  \\
  &=&   1_R \ot \mu^{-1}(  p_{2,\Z}^{\gamma}(a^m) - {m  \choose 2} p_{2,\Z}^{\gamma}(a) \,
p_{2,\Z}^{\gamma}(a) ) \\
  &=&   1_R \ot {m  \choose 2} a \widehat{\ot} a\,.
\end{eqnarray*}
If $m$ is odd, ${m  \choose 2}$ is divisible by $m$, so $\tau_3(a)=0$. Now suppose  that $m$ is
even. If $a=xy$, $y\in A_{(m/2)}$, $x \in A_{(m)}$ such that $x=\tilde{x}^2$ for some $\tilde{x}\in
A$, then 
  \begin{eqnarray*}
\tau_3(a)  &=&  1_R \ot 2{m  \choose 2} x \widehat{\ot} \tilde{x}  +  
1_R \ot 2{m  \choose 2} x \widehat{\ot} y  +  1_R \ot (m-1) \frac{m}{2} y \widehat{\ot} y \\
 &=& 0\,.
\end{eqnarray*}
To prove the converse inclusion, $\Ker{\tau_3} \subset ( A_{(m)} \cap A^2) A_{(m/2)}$, we may by a
standard argument suppose that $A$ is finitely generated. Choose a decomposition $t(A) =
\bigoplus_{p,r} \Z/p^r\Z \cdot a_{p,r}$ of the torsion subgroup $t(A)$ of $A$. Then any $a \in 
A_{(m)}$ can be written in the form $a = x^2 \prod_r a_{2,r}^{c_r}$ with $x\in A$ such that $x^2 \in 
A_{(m)}$, and $c_r \in \Z$ such that $c_r$ is odd if it is non-zero. As above, we get
    \begin{eqnarray*}
\tau_3(a)  &=&  \sum_r 1_R \ot {m \choose 2} c_r^2 \, a_{2,r} \widehat{\ot} a_{2,r}
+   \sum_{r<s} 1_R \ot 2{m \choose 2} c_r c_s  \, a_{2,r} \widehat{\ot} a_{2,s}  \\
 &=&  \sum_r 1_R \ot {m \choose 2} c_r^2 \, a_{2,r} \widehat{\ot} a_{2,r}\,.
\end{eqnarray*}
Now suppose  $a \in \Ker{\tau_3}$. Then    it follows from the decomposition SP$_2(A) =
\bigoplus_{p,r\le s} \Z/p^r\Z \cdot a_{p,r} \widehat{\ot} a_{p,s}$ that ${m \choose 2} c_r^2 \equiv 0$ 
mod $(2^r,m)$ for all $r$. This implies that $\frac{m}{2} c_r \equiv 0$ mod $2^r$ for all $r$, whence
$(\prod_r a_{2,r}^{c_r}) ^{\frac{m}{2}} = 1$. Thus $a \in ( A_{(m)} \cap A^2)
A_{(m/2)}$, as asserted. \hfbox \V

For $A=G/KN_2G^m$ we have $A_{(m )} = A$. Write $\bar{g} = g KN_2G^m$ for $g \in G$. Let  $a \in 
\Ker{\tau_3}$. Then by the lemma, $a = \bar{x}^2 \bar{y}$ where $x,y \in G$ such that $y=1$ if $m$
is odd, and $y^{\frac{m}{2}} \in KN_2G^m$ if $m$ is even. One has
  \[  \rho_{\ast}( 2 p_{2,\Z}^{\gamma}(\bar{x})  +  p_{2,\Z}^{\gamma}(\bar{y}) )  = a \]
and in $\Pz(G,KN_2G^m)$,
  \begin{eqnarray*}
m\, ( 2\kpz  ( {x})  +  \kpz  ( {y}) ) &=& \kpz(x^{2m}) -  {2m \choose  2} \bar{\mu}(\bar{x} \ot
xN_2)  \\
  & & \mbox{}  + 
\kpz(y^{ m}) -  { m \choose  2} \bar{\mu}(\bar{y} \ot yN_2)  \\
  &=&   \kpz(x^{2m}) + \kpz(y^{ m})  \\
  &=&   \kpz i(x^{2m}   y^{ m}) 
\,.
\end{eqnarray*}
Whence
  \[ 2p_{2,\Z}^{\gamma}(\bar{x})  +  p_{2,\Z}^{\gamma}(\bar{y})  =  P_{2,\Z}^{\cal N}(\pi)( 2\kpz 
( {x})  +  \kpz  ( {y})) \in  ( P_{2,\Z}^{\gamma}(G/KN_2G^m) )_{(m)}  \cap \rho_{\ast}^{-1}
\{a\}\,,\] and  
  \begin{eqnarray}\label{taueven} \bar{\tau}_2( a) = \tau_2 (2p_{2,\Z}^{\gamma}(\bar{x})  +  
p_{2,\Z}^{\gamma}(\bar{y}) ) &=&  1_R \ot (x^{2m}y^m U_m N_3)  \\
&=&   2\cdot  1_R \ot (x^m U_m N_3) +  1_R \ot (y^m  U_m N_3 ) \,. \nonumber
\end{eqnarray}
The latter equation    shows that for odd $m$, $\Imm{\bar{\tau}_2} = 1_R \ot (G^m U_m N_3 /U_m
N_3) $. Now abbreviate $V^{(m)} = G^m$ if $m$ is odd, and $V^{(m)} = G^{2m}V^m$ if $m$ is even.
Then  by   \REF{Imtau} and \REF{taueven} we obtain
  \[ \Ker{R\ot \kpz i} = \Imm{\tau_2}  = \Imm{\bar{\tau}_2}  = 1_R \ot (V^{(m)}U_mN_3 /U_m
N_3 )\,.\] Then the factorization
\comment{\\
\N\makebox[14.7cm]{ \makebox[0mm]{
\begin{minipage}{20cm}\small
  \[ p_{2,R}^{\cal N} i\,\colon\,  KN_2G^m \Sur{} R \ot  (KN_2G^m/N_3U_m) \mr{R  \ot p_{2,\Z}^{\cal N} i} R \ot
\Pz(G,KN_2G^m) \:\cong\: P_{2,R}^{\cal N}(G,KN_2G^m)  \]
\end{minipage}\ruled
}\rule[-5mm]{0mm}{3mm} }
\rule{0mm}{7mm}
}
\begin{eqnarray*}
\kpz i\,\colon\,  KN_2G^m  \Sur{}   R \ot  (KN_2G^m/U_mN_3) &\mr{R  \ot p_{2,\Z}^{\cal N} i}&
 R \ot \Pz(G,KN_2G^m) \\
&\cong&  P_{2,R}^{\cal N}(G,KN_2G^m)
\end{eqnarray*}
shows that 
  \[ KN_2G^m \cap (1 +  I_{R}(KN_2G^m)I_{R}(G) + I_{R,\cal N}^3(G)) =
\Ker{p_{2,R}^{\cal N} i} = U_m N_3 V^{(m)} \,,\]
also noting that $V^{(m)}$ contains  $(KN_2G^m)^m$. Together with \REF{Krepl} this proves the theorem. \hfbox

\comment{
Then the result
follows from commutativity of the following square. 
  \[ \matrix{
(G/N_2) \ot (KN_2/N_2)  & \mr{\nu}  & (G/N_2) \sm (G/N_2) \rule[-3mm]{0mm}{3mm} \cr
\mapdown{[\,,\,]}  &  & \mapdown{c}  \cr
((KN_3)\cap N_2)/N_3  &  \hra  & N_2/N_3 \rule{0mm}{7mm}
}\]
\Ph \hfill $\Box$\V
}

\section{The third relative Fox subgroup}\vspace{5mm}

We consider the following common
generalization of Fox subgroups and of relative
dimension subgroups:

\begin{defi}\label{foxdimdef}\quad   Let $G$ be a group and $H,K$ be subgroups of
$G$. For a commutative ring $R$ with unit $1_R$ let $I_R(G)$ denote the augmentation ideal of the group
algebra $R(G)$. Then define the {\bf `$n$-th relative Fox subgroup with respect to
$H,K,R$'}\/ to be the term
  \BE\label{Dndef}  G \cap (1_R + R(G) I_R(K) I_R(H) + I^{n}_R(G) I_R(H) )\:, \EE
with $n\ge 0$
and $I^0_R(G) = R(G)$. 
\end{defi}

Note that for $H=G$ and $K$ normal in $G$ this group is the $n+1$-st relative dimension
subgroup with respect to $K$  introduced by Passi (cf.\ \cite{Pa}). On the other hand, for $K=\{1\}$ we rediscover the classical $n$-th Fox subgroup of $G$ with respect to $H$. The mixed case ($H\neq G$ and $K\neq \{1\}$) seems to have been first studied by Tahara, Vermani and Razdan \cite{Razdan}  where the group \REF{Dndef} is determined for $n=2$, $R=\Z$, $H$ normal and a specific subgroup $K$ of $[H,G]$.\V

We determine the group  \REF{Dndef} for $n\le
2$ in full generality, as follows. The case $n=0$ is elementary; here
  \BE\label{n=1}  G \cap (1_R +  R(G) I_R(H)) = H \,, \EE
cf. \cite[Lemma 7]{Unico}.  For $n=1 $ we have the following.

\begin{prop}\label{D2F1}\quad Let $n_R$ denote the characteristic of $R$. Then 
the following is true, where we use the notation
of theorem \ref{D3R}.
\begin{Ar} 
\item If $n_R=0$ then
  \[  G \cap (1_R +  I_R(G) I_R(H)) = H_2  \prod_{p  \in \sigma(R)} p^et_p(\mbox{$H$ mod
$H_2$}) \,. \]
\item If $n_R > 0$ then
  \[  G \cap (1_R +  I_R(G) I_R(H)) = H_2 H^{n_R}\,.\]
\end{Ar}
\end{prop}

\proof Consider   the following sequence of homomorphisms
  \BE\label{n=2}
H/H_2   \to \frac{\dst R(G) I_R(H)}{\dst I_R(G) I_R(H)} \:\cong\: 
\frac{\dst   I_R(H)}{\dst  I_R^2(H)} \:\cong\: 
\frac{\dst  R \ot I_{\Z}(H)}{\dst \Imm{R \ot I_{\Z}^2(H)} }
\:\cong\:  R \ot \left( \frac{\dst I_{\Z}(H)}{\dst I_{\Z}^2(H)} \right)  \:\cong\:
  R \ot (H/H_2)   \EE
where the first one is given by $hH_2 \mapsto  h -1 +  I_R(G) I_R(H)$. The composition is the
canonical  morphism \map{j_R}{H/H_2}{R \ot (H/H_2)}. By   \REF{n=1}  we have 
$ (G \cap (1_R +  I_R(G) I_R(H)))/H_2 = \Ker{j_R} $, which was computed in \cite[Lemma 6]{Unico}. The
formula provided there gives the result.  \hfbox\V

For $n=2$ one uses the universal coefficient decomposition
obtained in \cite[Corollary 3]{Unico} to reduce from an arbitrary coefficient ring $R$ to the case
that $R=\Z/m\Z$, $m\ge 0$. Then the computation is completed by the following
result.

\begin{satz}\label{D3F2}\quad Let $G$ be a group and $H,K$ be subgroups of $G$. Let
$m\ge 0$ and $R=\Z/m\Z$. Then

\begin{Rom}

\comment{
\item \Ph\hfill $G \cap (1_R  + I_R(G) I_R(H) ) = [H,H] H^m
\:.$\hfill\Ph  
}

\item \Ph\hfill  $ G \cap (1_R + I_R(K) I_R(H) + I^2_R(G) I_R(H) ) = S_m\:,$\hfill\Ph
  \[ S_m
\:\stackrel{def}{=}\: {\rm sgp}\{ \prod_{h,k \in H} [h,k]^{a_{hk}} g^m \:\|\: g= \prod_{l\in
H}l^{b_l}\,,  \mbox{ all $ a_{hk}, b_l \in \Z$, and $ \forall k \in H \colon  $ }  \] 
  \[
\mbox{$\exists  d_k \ge 0 \,\colon \:k^{d_k} \in H_2H^m$  and } \prod_{h \in H} h^{a_{hk} - a_{kh} + 
{m \choose 2} b_h b_k}  \in K G_2 G^{d_k} \} \:. \]

\item If $H/H_2 $ is finitely generated we have the
following improvement of (i). Choose a decomposition $H/H_2  H^m \:\cong\: \bigoplus_{k=1}^r
\Z/d_k\Z \cdot(h_k H_2 H^m)$, $h_k \in H$. Then
  \[ G \cap (1_R + I_R(K) I_R(H) + I^2_R(G) I_R(H) ) = S_m^{fg}
H_3 H^{m^2}\:,\]
   \[ S_m^{fg} \:\stackrel{def}{=}\: {\rm sgp}\{ \prod_{1\le i<j\le
r}[h_i,h_j]^{a_{ij}}  (\prod_{l=1}^r h_l^{b_l})^m \:\|\:  a_{ij}, b_l \in \Z,\; \forall \:1 \le k\le
r \colon \]
  \[ h_k^{ {m \choose 2}b_k^2} \prod_{i<k} h_i^{a_{ik} + {m \choose 2} b_i b_k} \prod_{j>k}h_j^{\mbox{}
- a_{kj} + {m \choose 2} b_j b_k}    \in KG_2 G^{d_k} \} \:. \]

\end{Rom}
\end{satz}

\begin{bem}\label{Scontains}\quad \rm (1) It is easy to check directly that $ G \cap (1_R + I_R(K) I_R(H) + I^2_R(G) I_R(H) )$
contains
the canonical subgroup $H_3   V^{(m)}_H \,T_1T_2$ where $V^{(m)}_H = H^m$ if $m$ is odd and  $V^{(m)}_H = H^{2m}W^m$ with  $W=
\{h\in H
\,\|\,h^{m/2} \in KG_2G^m \}$ if $m$ is even, 
   \[ T_1 = {\rm sgp}\{\, [h,k^q]\:\|\,\mbox{\rm $h,k \in H$, $q\in \Z$, $h^q,k^q \in
KG_2G^m$}\}\:,\]
  \[  T_2 = {\rm sgp}\{\, [h,k]\:\|\,\mbox{\rm $h \in H\cap KG_2G^m$, $k \in H\cap KG_2G^mG^q$,
$q\in \Z$, $h^q  \in H_2$}\}\:.\]
 Note that $T_1$ contains $[H\cap KG_2G^m\,,H\cap KG_2G^m]$.\V

\N(2) Compared to the description of $ G \cap (1_{Z\!\!\!Z}  + I^2_{Z\!\!\!Z}(G) I_{Z\!\!\!Z}(H) )$ in \cite{Gu-Cu} for finitely generated $G$ the one given in
theorem \ref{D3F2}(ii) - apart from being more general (only $H/H_2H^m$ finitely
generated instead of $G/G_2$, arbitrary $R$ and $K$) - has the advantage not to require the choice of elementary-divisor-{\em compatible generators}\/ of $G/G_2$ and
 $HG_2/G_2$ but just the choice of {\em any basis}\/ of $H/H_2H^m$.

\end{bem}

\comment{
\begin{bem}\label{Scontains}\quad $S_m$ and $S_m^{fg}H_3H^{m^2}$ contain the canonical subgroup
  \[ H_3   V^{(m)}_H \, {\rm sgp}\{\, [h,k^q]\:\|\,\mbox{\rm $h,k \in H$, $q\in \Z$, $h^q,k^q \in
KG_2G^m$}\}\,,\]
 where  $V^{(m)}_H = H^m$ if $m$ is odd and  $V^{(m)}_H = H^{2m}W^m$,  $W= \{h\in H
\,\|\,h^{m/2} \in KG_2G^m \}$ if $m$ is even. In particular, $S_m  \supset H_3  (H_2)^m$, so one
does not need to care about the order of   the factors in   the products $\prod_{h,k \in H}
[h,k]^{a_{hk}}$ and $\prod_{l\in H}l^{b_l}$ (to see this in the latter case, also use the
Hall-Petresco formula $(l_1l_2)^m \equiv l_1^m l_2^m[l_2,l_1]^{{m \choose 2}}$   mod  $H_3$).
\end{bem}

\proof First consider $S_m$. Let $h_0\in H$ and $k_0 \in H_2$ or $h_0,k_0 \in H$, $q \in
\Z$ such that $h_0^q, k_0^q \in KG_2G^m$. In order to show that  $[h_0,k_0]$ or $[h_0,k_0^q]  \in
S_m$,  let $a_{h_0k_0}=1$ or $a_{h_0k_0}=q$, respectively, and let $a_{hk} = 0$ for $(h,k) \neq
(h_0,k_0)$  and   all $b_l = 0$.    Then the two conditions in the definition of $S_m$ are
satisfied since  \begin{itemize}
\item in the case $h_0  \in H$ and $k_0 \in H_2$,   one can take $d_{k_0} =1$ and $d_k = 0$ for
$k\neq 0$;

\item in the case $h_0^q,k_0^q \in KG_2G^m$, one can take $d_k = m$ for all $k\in H$; it remains
to observe that $[h_0,k_0]^q \equiv [h_0,k_0^q]$ mod $H_3$. \end{itemize}

Now let $g\in H$ and suppose in addition that $g^{m/2}\in KG_2G^m$ if $m$ is even. To show that $g^m 
\in S_m$, let $b_g =1$ and $b_l =0$    for $l\neq g$, and all $a_{hk}=0$. Then taking $d_k =m$ for all
$k\in H$ the conditions are satisfied since $m$ divides ${m\choose 2}$ if $m$ is odd, and since
$g^{{m \choose 2}} = (g^{m/2})^{m-1} \in KG_2G^m$  otherwise.

Suppose $m$ is even. Let $c\in H$. To show that $c^{2m} \in  S_m$     let $b_c =2$ and $b_l  =0$
otherwise, and all $a_{hk}=0$. Then    taking $d_k = m$  for all $k\in H$ does since $m$ divides $
{m \choose 2}2^2$.

Now consider $S_m^{fg}$. The relation $S_m^{fg}H_3H^{m^2} \supset (H \cap KG_2G^m)_2$ is a consequence  of
the theorem  since $G \cap (1_R + I_R(K) I_R(H) + I^2_R(G) I_R(H) ) \supset (H \cap KG_2G^m)_2$   
which  follows   from commutativity  of diagram \REF{cHGdia}   below. To show that  
$S_m^{fg}H_3H^{m^2} \supset V^{(m)}_H$ we first check that $S_m^{fg}H_3H^{m^2} \supset (H_2)^m$.
Indeed, modulo $H_3$ any element of $(H_2)^m$ is of the form $\prod_{1\le i<j\le
r}[h_i,h_j]^{a_{ij}}    $ where $m$ divides all exponents $a_{ij}$. Thus taking all $b_l=0$ we have  
  \[ \prod_{i<k} h_i^{a_{ik}}   \prod_{j>k}h_j^{\mbox{} - a_{kj}}    \in
H_2H^m  \subset G_2G^{d_k}  \]
   for all $k$ since $d_k$ divides $m$. Whence $S_m^{fg}H_3H^{m^2} \supset
(H_2)^m$. Now let $h\in H$  such that, if $m$ is even, $h^{m/2}\in KG_2G^m$. Modulo $H_2H^m$, $h$ is
of the form $\prod_{l=1}^r h_l^{b_l}$, whence $h^m \equiv (\prod_{l=1}^r h_l^{b_l})^m$  mod 
$H_3H^{m^2}(H_2)^m$. So for  proving $h^m \in S_m^{fg}H_3H^{m^2}$ it suffices to show  that
$(\prod_{l=1}^r h_l^{b_l})^m \in S_m^{fg}$. Taking all $a_{ij}=0$    we have for all $k$:
  \begin{eqnarray}
 \pi_k &\stackrel{def}{=}  & h_k^{ {m \choose 2}  b_k^2} \prod_{i<k} h_i^{  {m \choose 2} b_i b_k}
\prod_{j>k}h_j^{ 
  {m \choose 2} b_j b_k}      \nonumber \\
 &\equiv&  (\prod_{l=1}^r h_l^{b_l})^{{m \choose 2} b_k} \quad \mbox{mod $H_2$} \nonumber \\
   &\equiv& h^{{m \choose 2} b_k}  \quad \mbox{mod $H_2H^m \subset G_2 G^{d_k}$} \,.
\label{auskl}
\end{eqnarray}
But $h^{{m \choose 2} b_k}$ is in  $G^{m} \subset G^{d_k}$ if $m$ is odd, and in $KG_2G^m \subset
KG_2G^{d_k}$ if $m$ is even, since then $ h^{{m \choose 2} b_k}  = (h^{m/2})^{(m-1)b_k}$.

Finally, suppose  $m$ is even,   and let $h\in H$. Write $h    \equiv \prod_{l=1}^r h_l^{c_l}$ mod
$H_2H^m$; then $h^2 
\equiv \prod_{l=1}^r h_l^{2c_l}$. 
In order to prove that $h^{2m} \in S_m^{fg}H_3H^{m^2}$
it suffices   to show that $ (\prod_{l=1}^r h_l^{2c_l} )^m \in S_m^{fg}$. So let $b_l=2c_l$,   all
$a_{ij}=0$.  Then by \REF{auskl}, $\pi_k \equiv (h^m)^{(m-1)c_k} \equiv  0 $   mod  $ G_2G^m \subset
G_2G^{d_k}$   for all $k$, which proves the assertion.   \hfill $\Box$\vspace{4mm}
}

\proofofthm{D3F2}
In the sequel we shall frequently use the remarks around \REF{XotR}.
Consider the following commutative diagram, where we abbreviate $(H/H_2)\hsm{2} =
(H/H_2)\sm (H/H_2)$.
 \BE\label{cHGdia} \matrix{
\ruleu\ruled 
 (H/H_2)\hsm{2}  &  \mr{\ell} &   (H/H_2)\ot (H/H_2) \ot R  & \mr{i\ot id \ot R}  &  (G/KG_2)
\ot (H/H_2) \ot R \cr 
\surdown{c}  &  & \surdown{\mu_{H,R}} &  &  \surdown{\mu_{G,R}}\cr
\rule{0mm}{6mm} H_2/H_3  &  \mr{p_{2,R}}  &  I_R^2(H) / I_R^3(H)  &  \mr{j}  & 
\frac{\dst I_R(G) I_R(H)}{\dst  I_R(K)I_R(H) + I^2_R(G)   I_R(H) }\rule[-5mm]{0mm}{3mm} 
}\rule[-5mm]{0mm}{3mm}\EE 
\rule{0mm}{7mm} where for $a,b,h \in H$, $g \in G$,   and $r\in R$
  \[ l((aH_2) \sm (bH_2)) = (aG_2) \ot (bH_2) \ot 1_R - (bG_2) \ot (aH_2) \ot 1_R \:,\]
  \[
 c((aH_2) \sm (bH_2)) = aba^{-1}b^{-1}H_3\:,\]
\[  
\mu_{H,R}((aG_2)\ot (bH_2) \ot r) = r(a-1)(b-1) + I_R^3(H) \:, \]
  \[ \mu_{G,R}((gG_2)\ot (hH_2) \ot r) = r(g-1)(h-1) + 
I_R(K)I_R(H) + I^2_R(G) I_R(H)\:, \]
and where $i$ and $j$ are induced by the inclusion $H \hra G$. 

By  proposition \ref{D2F1}, $G \cap (1_R +  I_R(G) I_R(H)) = H_2H^m$. Modulo $I_R(K) I_R(H) +
I_R^2(G) I_R(H)$,  we have for $h\st \in H_2, h \in H$
  \begin{eqnarray}
h\st h^m - 1_R  &\equiv&  (h\st - 1_R) + m (h  - 1_R) + { m \choose 2 } (h  - 1_R)^2 \nonumber\\
 &=&\label{p2h'h}  j  \Big( p_{2,R} (h\st) + { m \choose 2 } \mu_{H,R}(hH_2 \ot hH_2 \ot 1_R)
\Big)\,. \end{eqnarray}
Now we use the crucial fact (\cite{Fox2}, proof of theorem 3.1)  that the right-hand square of
diagram \REF{cHGdia} is a cocartesian (or pushout) square which implies that $\Ker{j}= \mu_{H,R}\Ker{i\ot id \ot R}$.
  Using this and
 proposition \ref{D2F1} we conclude that $x \in G \cap (1_R +
I_R(K)I_R(H) + I_R^2(G)I_R(H) )$ {\em if and only if}\/ $x = h\st h^m$ as above such that in
$I_R^2(H) / I_R^3(H)$,
  \BE \label{p2R=mu}  p_{2,R} (h\st h^m)  = p_{2,R} (h\st) + { m \choose 2 } \mu_{H,R}(hH_2 \ot
hH_2\ot 1_R) \,\in\, \mu_{H,R}\,\Ker{i\ot id \ot R}\,.\EE
This situation is further analyzed in the following two lemmas.

\begin{lem}\label{Y}\quad The preimage of   the subset $\kpzr(H_2H^m) \subset I^2_R(H)/I^3_R(H)$
under the map $\mu_{H,R}$ is 
  \[  Y  \stackrel{def}{=} \Imm{\ell} + \{ {m\choose 2} (hH_2  \ot  hH_2  \ot 1_R)\,\|\, h\in H\}
\:\subset\: (H/H_2)  \ot (H/H_2) \ot R \,.\]
\end{lem}

\proof  First of all, $\mu_{H,R}(Y) = \kpzr(H_2H^m)$; this follows  from \REF{p2h'h} and the diagram
 above. Secondly, we show that $Y$ contains $\Ker{\mu_{H,R}}$. For this purpose, recall
the isomorphism $I^2_{\Z}(H)/I^3_{\Z}(H) \cong {\rm U}_2{\rm L}(G) = ((H/H_2)  \ot
(H/H_2))/\ell\,\Ker{c}$ from   \cite{Ba-GrII}. Using the six-term exact sequence  for the tensor
product  one deduces exact sequences
  \[ \Ker{c}   \mr{\ell} (H/H_2)  \ot (H/H_2) \ot R  \mr{\mu_{H,\Z} \ot R} 
(I^2_{\Z}(H)/I^3_{\Z}(H))  \ot R   \to 0 \]
  \[  {\rm Tor}_1^{\Z}(H/H_2, R)  \mr{\tau} (I^2_{\Z}(H)/I^3_{\Z}(H))  \ot R \mr{\nu \ot R}
(I_{\Z}(H)/I^3_{\Z}(H))  \ot R \mr{\pi \ot R} (H/H_2)  \ot R \to 0 \,,\]
where $\nu$ is the inclusion and $\pi$ is the quotient map modulo $\Imm{\nu}$, combined with  the
canonical isomorphism $I_{\Z}(H)/I^2_{\Z}(H) \,\cong\, H/H_2$. 
  Identifying ${\rm Tor}_1^{\Z}(H/H_2, R)$ with the set of $m$-torsion elements of $H/H_2$,
let $h\in H$  such that $h^m \in H_2$. Then 
  \[ \tau(hH_2)  = \nu^{-1}(m p_{2,\Z}(h) ) \ot 1_R = \Big(p_{2,\Z}(h^m) -  {m \choose 2}
\mu_{H,\Z}(hH_2  \ot hkH_2) \Big) \ot 1_R   \,.\]
Using the canonical isomorphism $(I_{\Z}(H)/I^3_{\Z}(H))  \ot R \,\cong\, (I_{R}(H)/I^3_{R}(H)) $,
we deduce that $\Ker{\mu_{H,R}} = \Ker{(\nu \ot R)(\mu_{H,\Z} \ot R)} \subset Y$.

It remains to show that $Y$ is a subgroup.  Let $x,y \in (H/H_2)  \sm (H/H_2)$  and $h,k \in H$. Then
  \begin{eqnarray*}
& & \Big( \ell(x)  + {m \choose 2} (hH_2 \ot hH_2 \ot 1_R) \Big)   -  \Big( \ell(y)  +  {m \choose 2} (kH_2
\ot kH_2 \ot 1_R ) \Big) \\
&=&   \ell(x-y) +  {m \choose 2}( hkH_2 \ot hkH_2 \ot 1_R ) +  {m \choose 2}
(hH_2 \ot kH_2 \ot 1_R) \\
  & & - {m \choose 2} ( kH_2 \ot hH_2 \ot 1_R)  - 2 {m \choose 2} ( kH_2 \ot kH_2 \ot 1_R)   \\
 &=&  \ell\Big(x-y  +  {m \choose 2} ( h H_2 \sm  kH_2)\Big) +   {m \choose 2} ( hkH_2 \ot hkH_2 \ot
1_R)   \\
 &\in& Y\,.
\end{eqnarray*}

Thus the lemma is proved.\hfill $\Box$\V

\begin{lem}\label{T}\quad One has $Y\cap \Ker{i\ot id \ot R} = T_m$, with
\comment{
  \[  T_m = \Big\{  \ell(\sum_{h,k \in H} {a_{hk}}  (hH_2) \sm (kH_2))   + {m \choose 2} (gH_2) \ot
(gH_2) \ot 1_R \:\|\: g= \prod_{l\in
H}l^{b_l}\,,  \mbox{ all $ a_{hk}, b_l \in \Z$,   }  \] 
  \[
\mbox{  and $ \forall k \in H \colon \exists  d_k \ge 0 \,\colon \:k^{d_k} \in H_2H^m$  and }
\prod_{h \in H} h^{a_{hk} - a_{kh} +  {m \choose 2} b_h b_k}  \in KG_2 G^{d_k} \Big\} \:. \]
}
 \[  T_m = \Big\{ \ell \Big(\sum_{h,k \in H} {a_{hk}}  ( hH_2  \sm    kH_2 )  \Big)  + {m \choose 2}
(gH_2  \ot
 gH_2  \ot 1_R ) \:\|\: g= \prod_{l\in
H}l^{b_l}\,,  \mbox{ all $ a_{hk}, b_l \in \Z$,   }  \] 
  \[
\mbox{  such that $ \forall k \in H \colon \exists  d_k \ge 0 \,\colon \:k^{d_k} \in H_2H^m$  and }
\prod_{h \in H} h^{a_{hk} - a_{kh} +  {m \choose 2} b_h b_k}  \in KG_2 G^{d_k} \Big\} \:. \]
If $H$ is finitely generated, we have (with the notation of \ref{D3F2} (ii)) 
  \[  Y\cap \Ker{i\ot id \ot R} = T_m^{fg} \,, \]
  \[ T_m^{fg} =  \{ \ell \Big( \sum_{1\le i<j\le
r} {a_{ij}}  (h_iH_2  \sm  h_jH_2) \Big)  +  {m \choose 2} (gH_2  \ot  gH_2  \ot 1_R )
  \:\|\: g = \sum_{l=1}^r   h_l^{b_l} \mbox{ such that }  \]
  \[ \forall \:1 \le k\le r\,
\colon\, h_k^{ {m \choose 2}  b_k^2} \prod_{i<k} h_i^{a_{ik} + {m \choose 2} b_i b_k}
\prod_{j>k}h_j^{\mbox{} - a_{kj} + {m \choose 2} b_j b_k}  \in KG_2 G^{d_k} \} \:. \]
\end{lem}

\proof Let $y \in Y$, i.e. $y = \ell(\sum_{h,k \in H} {a_{hk}}       (hH_2 \sm  kH_2))   + 
{m \choose 2} (gH_2   \ot  gH_2  \ot 1_R)$ for some $a_{hk} \in \Z$, $g \in H$. Consider an
arbitrary decomposition $g = \prod_{l\in H}l^{b_l}$, $b_l \in \Z$. Using the isomorphism
$\alpha\,\colon\,(G/KG_2)  \ot (H/H_2) \ot R \,\cong\, (G/KG_2)  \ot (H/H_2H^m) $ we obtain
     \begin{eqnarray}
\alpha(i\ot  id \ot R)(y)  &=& \sum_{h,k \in H} {a_{hk}} \Big( (hKG_2) \ot (kH_2H^m)  - (kKG_2) \ot
(hH_2H^m) \Big) \nonumber\\
  & &   +  {m \choose 2} \sum_{l_1,l_2 \in H} b_{l_1} b_{l_2}  (l_1KG_2) \ot (l_2H_2H^m)  \nonumber\\
  &=&  \label{zerleg} \sum_{ k \in H} \Big(\prod_{h \in H} h^{a_{hk} - a_{kh} + 
{m \choose 2} b_h b_k} KG_2 \Big) \ot (k H_2 H^m)
\end{eqnarray}

So if $y$ and the chosen decomposition of $g$ satisfy the condition in the definition of $T_m $ then
it follows from bilinearity of the tensor product that $(i\ot  id \ot R)(y) = 0$. Whence $T_m \subset
Y\cap \Ker{i\ot id \ot R}$. To prove the converse inclusion we may by a standard argument suppose
that $H$ is finitely generated. Then $y \in Y$ can be written in the form
  \[ y = \ell\Big( \sum_{1\le i<j\le
r} {a_{ij}}  (h_iH_2  \sm  h_jH_2) \Big)  +  {m \choose 2} (gH_2  \ot  gH_2 \ot 1_R) \]
with $g = \prod_{l=1}^r   h_l^{b_l}$, all $a_{ij}, b_l \in \Z$. To see this just observe that $
 ( H/H_2) \ot  (H/H_2) \ot R \:\cong\: ( H/H_2H^m) \ot  (H/H_2H^m)$ and that $\ell$ factors through
$( H/H_2H^m) \sm  (H/H_2H^m)$. Now we use the composite isomorphism   
  \[ \Psi\,\colon\,  ( G/KG_2) \ot  (H/H_2) \ot R \:\stackrel{\alpha}{\cong} \: ( G/KG_2) \ot 
(H/H_2H^m) \:\cong\: \]
  \[  \bigoplus_{k=1}^r  ( G/KG_2) \ot (\Z/d_k\Z \cdot(h_kH_2H^m) ) \:\cong\:  \bigoplus_{k=1}^r    
G/(KG_2 G^{d_k}) \] By \REF{zerleg} we have 
 \[  \Psi(i\ot  id \ot R)(y) =  \bigoplus_{k=1}^r  \left( h_k^{ {m \choose 2}  b_k^2} \prod_{i<k}
h_i^{a_{ik} + {m \choose 2} b_i b_k} \prod_{j>k}h_j^{\mbox{} - a_{kj} + {m \choose 2} b_j b_k} \,
KG_2 G^{d_k} \right)  \,.\]
So $Y\cap \Ker{i\ot id \ot R} = T_m^{fg}$. Moreover, $T_m^{fg} \subset T_m$;  to see this,  let
$a_{hk} = b_l = 0$ unless $h,k,l \in \{h_1,\ldots,h_r\}$, and $a_{h_ih_j} = 0$ if $i>j$. Then
an element of $T_m^{fg}$ satisfies the conditions defining  $T_m$ in taking   $d_k=m$ if
$k\notin \{h_1,\ldots,h_r\}$ and $d_{h_k}  =  d_k$   for $1 \le k\le r$. Thus also the inclusion $
T_m \subset Y\cap \Ker{i\ot id \ot R} $ already proved is an equality, as asserted. \hfbox\V

Summarizing the above we obtain equations
  \[ p_{2,R}(G \cap (1_R + I_R(K) I_R(H) + I^2_R(G) I_R(H) )) = \mu_{H,R}(T_m) =  p_{2,R}(S_m) \,,
\] analogously for $T_m^{fg}$ and $S_m^{fg}$, where the first equation follows from \REF{p2R=mu},
Lemmas \ref{Y} and \ref{T}, and the second one from \REF{cHGdia} and \REF{p2h'h}. So it
only remains to observe that $\Ker{\map{p_{2,R}}{H_2H^m}{P_{2,R}(H)}}  = D_{3,R}(H)$ is contained
in $V^{(m)}_H$ by Sandling's formula for $D_{3,R}(H)$, see \cite{Sa}  (or theorem \ref{D3} above
for $K=\{1\}$), and that $V^{(m)}_H$ is contained in $S_m$ and  $S_m^{fg}$ by remark
\ref{Scontains}. Thus theorem \ref{D3F2} is proved. \hfbox
\vspace{5mm}

\begin{center} \sc Acknowledgment  \end{center} \vspace{3mm}

{\small I am grateful to C.\ K.\ Gupta for sending me a copy of her preprint on the second Fox
subgroup which was helpful to overcome a certain technical difficulty at some point of my own research
on the topic.}

\vspace{2cm}

 
 \N LAMAV, ISTV2\\
Universit\'e de Valenciennes et du Hainaut-Cambr\'esis,\\
Le Mont Houy,\\
59313 Valenciennes Cedex 9\\
 France

\tw

\N U.R.A. au C.N.R.S. 751, et\\
D\'epartement de Math\'ematiques,\\
Universit\'e de Valenciennes, \\ 
le Mont Houy, BP 311,\\
59304 Valenciennes Cedex, France.\\
Email: Manfred.Hartl@univ-valenciennes.fr
\vspace{2mm}

\N email: hartl@univ-valenciennes.fr

\tw